\documentclass[final,3p,times]{elsarticle}

\usepackage{amssymb}
\usepackage{amsmath}
\usepackage{amsthm}

\usepackage{bm}
\usepackage{algorithm}
\usepackage{algpseudocode}
\usepackage{svg}

\usepackage{tikz}
\usepackage{booktabs}
\usepackage{subcaption}

\usepackage{hyperref}
\usepackage{url}
\usepackage{enumitem}

\begin{document}

\begin{frontmatter}



\title{Physics-constrained identification of graph-based thermal networks for spacecraft digital twins}

\author[1]{Luca Sosta}
\author[2]{Carlo Ciancarelli}
\author[2]{Leonardo Marini}
\author[1]{Stefano Pagani}
\author[1]{Francesco Regazzoni}
\author[1]{Nicola Parolini}

\affiliation[1]{organization={MOX-Dipartimento di Matematica, Politecnico di Milano},
               addressline={Piazza Leonardo da Vinci 32},
               city={Milan},
               country={Italy}}

\affiliation[2]{organization={Thales Alenia Space Italia},
               addressline={Via Saccomuro 24},
               city={Rome},
               country={Italy}}

\begin{abstract}
Reconstructing a thermal model capable of efficiently simulating the behavior of a spacecraft from sparse and localized temperature measurements remains a challenging task. To address this,
we introduce a physically-constrained calibration framework for Lumped Parameter Thermal Models (LPTMs), formulated as a trajectory-based inverse problem for graph dynamical systems. The model reconstructs thermal dynamics directly from temperature measurements and known inputs, without relying on a priori parameter values derived from material properties or geometric assumptions.
Physical admissibility is enforced at the parameterization level: positivity of nodal coefficients and symmetry of conductive interactions are imposed by construction. This guarantees stable dynamics and restricts the identification problem to a physically meaningful parameter space, improving conditioning without the need of additional regularization. 
The identification problem is addressed through trajectory matching, ensuring stable rollout over extended time horizons.
The methodology is validated on synthetic datasets generated from high-fidelity finite element simulations under progressively complex forcing conditions. The calibrated LPTMs accurately reproduce long-term temperature evolution and exhibit robustness to measurement noise. 
The proposed framework provides a systematic approach to the calibration of reduced-order thermal models by combining physical structure with data-driven identification. The numerical results show a favorable balance between accuracy and computational  efficiency, making the models suitable for integration in spacecraft thermal Digital Twin applications.
\end{abstract}




\end{frontmatter}



\section{Introduction}

Thermal modeling of spacecraft involves the prediction of energy transfer in heterogeneous environments governed by the heat equation and radiative exchange mechanisms \cite{carslaw1959conduction,incropera2007fundamentals,modest2013radiative}. In the absence of atmospheric convection, satellites are subject to strongly time-dependent boundary conditions, including solar irradiation, albedo, and eclipse phases, leading to nonlinear and non-stationary temperature dynamics. As a result, inadequate thermal management under such conditions may compromise subsystem functionality, reduce mission lifetime, or ultimately lead to system failure.

The Digital Twin (DT) paradigm, initially developed in the aerospace sector to support real-time monitoring and predictive maintenance of complex systems \cite{glaessgen2012digital}, aims to combine physics-based models with real-time data for monitoring, prediction, and control. In thermal applications, a DT must remain accurate over long prediction horizons while retaining computational efficiency compatible with real-time deployment and robustness across varying operating conditions.

\smallskip

High-fidelity models based on the finite element method (FEM) provide accurate solutions but are computationally expensive and unsuitable for real-time DT applications. Reduced-order modeling techniques offer a systematic alternative by approximating high-dimensional systems through low-dimensional representations while preserving essential dynamical properties. In this context, projection-based model order reduction techniques derived from FEM discretizations have been proposed to accelerate large-scale simulations while preserving physical fidelity \cite{antoulas2005approximation,benner2015survey,quarteroni2016rbm,hesthaven2015rbm}. Lumped Parameter Thermal Models (LPTMs) can be also interpreted as reduced-order models of the underlying heat equation, where the thermal field is approximated by a network of interacting nodes. From a structural perspective, LPTMs can be interpreted as diffusive dynamical systems defined on graphs, where thermal interactions are governed by Laplacian operators. Such systems belong to the broader class of network dynamical systems, extensively studied in control theory and multi-agent systems \cite{mesbahi2010graph}. However, unlike projection-based methods, their parameters usually do not directly inherit physical meaning (except in very simple settings, with homogeneous materials, simple geometries, and clearly defined heat transfer paths) and must be inferred from data \cite{torralbo2018thermal}.

Data-driven approaches, including neural networks and operator learning methods \cite{lu2021deeponet}, have been proposed to approximate complex dynamical systems at reduced computational cost. In recent years, physics-informed machine learning has expanded the range of possible approaches by embedding physical structure into learning architectures to enhance robustness and generalization in dynamical systems \cite{raissi2019pinn,karniadakis2021piml}. In parallel, system identification techniques for continuous-time dynamics have been revisited through sparse regression methods such as SINDy \cite{brunton2016discovering} and Neural Ordinary Differential Equations (Neural ODEs) \cite{chen2018neuralode}. Hybrid approaches, such as Universal Differential Equations (UDEs), further combine mechanistic models with data-driven components, enabling the integration of prior physical knowledge within the learning process \cite{rackauckas2020universal}. Structure-preserving approaches such as Hamiltonian neural networks have also been proposed to enforce physical invariants in learned models \cite{greydanus2019hnn}. Among these, reduced representations such as Dynamic Mode Decomposition (DMD) have been successfully applied to approximate complex dynamical systems at lower computational cost \cite{schmid2010dynamic,tu2014dynamic}. Despite their flexibility, these approaches often suffer from poor extrapolation capabilities, lack of physical interpretability, and instability over long prediction horizons, with stability preservation becoming particularly critical in data-scarce regimes \cite{reis2011lyapunov}.
\smallskip

These limitations highlight a gap between physically interpretable reduced-order models and flexible data-driven approaches. In particular, the identification of reduced thermal models from sparse measurements remains an ill-posed inverse problem, due to limited observability and non-uniqueness of the parameters \cite{stuart2010inverse,kaipio2005inverse}.

In this work, we propose a physics-constrained identification framework for lumped thermal models formulated as graph dynamical systems. The approach combines:
(i) a graph-based representation of thermal interactions,
(ii) a continuous-time system identification formulation inspired to Neural ODEs,
(iii) a structure-preserving parameterization enforcing positivity, symmetry, and dissipativity by construction.
This formulation ensures that the identified model remains within the class of physically admissible thermal systems, improving stability and generalization without requiring explicit regularization.
\smallskip

The remainder of this paper is organized as follows. Section \ref{sec:2} introduces the theoretical background, including the LPTM formulation, its graph representation, and the connection with Neural ODE models. Section \ref{sec:3} presents the inverse calibration framework and the associated learning problem. Section \ref{sec:4} describes the computational architecture and the enforcement of physical constraints. Section \ref{sec:5} reports numerical results on synthetic test cases, including sensitivity analyses and robustness assessments. Section \ref{sec:6} concludes the paper and outlines future research directions.
\section{Background and Problem Formulation}
\label{sec:2}

We consider a thermal subsystem equipped with a finite number of temperature probes placed at specific locations within the structure. 
These sensors provide temperature measurements over time and represent the primary information available for reconstructing the thermal dynamics of the system. Let $\{t_k\}_{k=0}^{N}$ denote the discrete sampling times. At each time instant, the probes provide a vector of temperature measurements $y_k \in \mathbb{R}^p$,
where $p$ is the number of probes.
Internal heat sources are assumed to act at the same probing points. 
These sources account for heater activations or thermal dissipation generated by electronic components, whose power levels are typically known from telemetry. 
Their values are collected in the input vector $u_k \in \mathbb{R}^p $ and the available dataset can be written as
$\mathcal{D} = \{u_k, y_k\}_{k=0}^{N}$.
External thermal interactions, such as radiative exchanges with the surrounding environment, are not directly included in this description and will be introduced separately in the following sections.

We distinguish between the number of measured outputs \(p\) and the number of thermal states \(n\) used in the model. In general, the graph-based thermal model may include additional nodes beyond sensor locations, corresponding to lumped components or internal thermal masses. Therefore, the state dimension \(n\) does not necessarily coincide with the number of measurements \(p\). The observation operator is implicitly assumed to map the full state \(T(t) \in \mathbb{R}^n\) to the measured outputs \(y(t) \in \mathbb{R}^p\). In the numerical experiments presented in this work, we consider the case \(n = p\), that is, one state per sensor, for simplicity.

\subsection{Lumped Parameter Thermal Models (LPTMs)}
\label{sec:lptm}

LPTMs are reduced-order models that approximate heat transfer problems by aggregating the domain into a finite set of control volumes (nodes) connected by lumped thermal couplings (edges). The evolution of the temperature $T$ in a domain $\Omega \subset \mathbb{R}^3$ is described by the heat equation arising from conservation of energy and Fourier's law of conduction \cite{carslaw1959conduction,incropera2007fundamentals}

$$
\rho c_p \frac{\partial T}{\partial t} = p + \nabla \cdot (\lambda \, \nabla T),
$$
where $\rho$ denotes the mass density, $c_p$ is the volume specific heat at constant pressure, $p$ represents any volumetric heat source and $\lambda$ is the conductivity. The domain $\Omega$ is partitioned into disjoint control volumes $\Omega_i$, such that $\Omega = \bigcup_{i=1}^n \Omega_i$. The node temperature is defined as the volume average over a control volume, namely

$$
T_i(t) := \frac{1}{|\Omega_i|} \int_{\Omega_i} T(\bm{x},t)\, dV,
$$

while the lumped thermal capacitance is defined as

$$
C_i := \int_{\Omega_i} \rho(\bm{x}) c_p(\bm{x})\, dV.
$$

Let $\Gamma_{ij} := \partial \Omega_i \cap \partial \Omega_j$ be the interface between nodes $i$ and $j$, with outward normal $\bm{n}_{ij}$ from $\Omega_i$ to $\Omega_j$. The net conductive heat flow from $i$ to $j$ is

$$
\Phi_{ij}(t) = \int_{\Gamma_{ij}} -\lambda(\bm{x})\, \nabla T(\bm{x},t)\!\cdot\! \bm{n}_{ij}\, dS.
$$

Under the usual lumped approximation, this interfacial flow is approximated as

$$
\Phi_{ij}(t) \approx G_{ij}\, \big(T_i(t)-T_j(t)\big),
$$

where $G_{ij}=1/R_{ij}>0$ is an effective conductance that aggregates geometry, materials, and contact properties. External heat inputs acting on node $i$ (heaters, absorbed radiation) are collected in $P_i(t)$. The nodal energy balance becomes

$$
C_i \, \frac{dT_i}{dt} = P_i\,(t) + \sum_{j\in \mathcal{N}(i)} G_{ij}\, \big(T_j - T_i\big), \qquad i=1,\dots,n,
$$

where $\mathcal{N}(i)$ is the set of neighbors of node $i$.

\subsection{Graph-based description: nodes, edges, and matrices}

We consider a thermal graph $G=(V,E)$, where the thermal system is represented as a network of interacting components \cite{chung1997spectral,newman2010networks}. This graph-based formulation is naturally connected to dynamical systems representations, including structural controllability and observability concepts \cite{lin1974structural}, and is closely related to recent graph-neural and message-passing approaches for multi-component physical systems. In this framework, nodes (vertices) represent control volumes $\Omega_i$ associated with physical components or subsystems. Each node is associated with a single temperature state $T_i$ and a capacitance $C_i$ corresponding to the integral of $\rho c_p$ over $\Omega_i$. Meanwhile, edges encode lumped thermal couplings between nodes. Each edge $(i,j)$ is weighted by a symmetric conductance $G_{ij}=G_{ji}>0$ that aggregates conduction paths. We define the weighted adjacency $W=[w_{ij}]$ with $w_{ij}=G_{ij}$ for $i\neq j$ and $w_{ii}=0$, the degree matrix $D=\mathrm{diag}(d_i)$ with $d_i=\sum_{j\neq i} G_{ij}$, and the (combinatorial) Laplacian $L$:
\begin{equation*}
L := D - W, \qquad L_{ij} = \begin{cases}
-\,G_{ij}, & i\neq j,\\
\sum_{k\neq i} G_{ik}, & i=j.
\end{cases}
\end{equation*}

For simplicity, we aggregate both controlled fluxes (heaters) and external excitations (orbital fluxes, disturbances) into a single input vector $u(t)$. Thus, we obtain:

\begin{equation}
C \, \dot T(t) = -L T(t) + u(t), \qquad C=\mathrm{diag}(C_1,\dots,C_n).
\label{eq:lptm}
\end{equation}

\subsection{Radiative fluxes and linearization}
In addition to conductive exchanges between nodes, satellites experience radiative heat transfer. Two main cases can be distinguished: 

\begin{itemize}[leftmargin=11pt]
    \item \textbf{Radiative exchange with an external sink (deep space)}\\
The net flux between a node $i$ at temperature $T_i$ and a radiative sink is governed by the Stefan--Boltzmann law \cite{modest2013radiative}:

$$
P_i^{\text{rad}} = \epsilon_i \sigma A_i \left(T_i^4 - T_{\infty}^4\right),
$$

where $\epsilon_i$ is the emissivity, $A_i$ the radiating area, and $\sigma$ the Stefan--Boltzmann constant. In most satellite applications, the sink is deep space with $T_{\infty} \approx 0$, so that

$$
P_i^{\text{rad}} \approx \epsilon_i \sigma A_i T_i^4.
$$

Typically, for satellite thermal monitoring, all quantities $A_i$, $\epsilon_i$, and $T_{\infty}$ are known or fixed by design, and the corresponding radiative conductance can be computed analytically without ambiguity. Therefore, this contribution does not require calibration. \\

\item \textbf{Radiative exchange between two components}\\
For two nodes $i$ and $j$ with view factor $F_{ij}$, emissivities $\epsilon_i, \epsilon_j$, and areas $A_i, A_j$, 
the net radiative power exchanged is:

\begin{equation*}
P_{ij}^{\text{rad}} =
\frac{\sigma\!\left(T_i^{4}-T_j^{4}\right)}
{\displaystyle \frac{1-\epsilon_i}{\epsilon_i A_i}\;+\;\frac{1}{A_i F_{ij}}\;+\;\frac{1-\epsilon_j}{\epsilon_j A_j}},
\end{equation*}

where $\sigma$ is the Stefan--Boltzmann constant.  
We introduce an effective emissivity $\epsilon_{ij}^\ast$ defined by:

\begin{equation*}
\frac{1}{\epsilon_{ij}^\ast A_i F_{ij}}
:=\frac{1-\epsilon_i}{\epsilon_i A_i}+\frac{1}{A_i F_{ij}}+\frac{1-\epsilon_j}{\epsilon_j A_j},
\end{equation*}

so that the net exchange can be compactly written as:
\begin{equation*}
P_{ij}^{\text{rad}}
= \sigma\,\epsilon_{ij}^\ast A_i F_{ij}\,\left(T_i^{4}-T_j^{4}\right),
\end{equation*}

An equivalent conductance between nodes $i$ and $j$, 
can be obtained linearizing the radiative exchange around a reference temperature $\bar T$ (e.g. the mean of $T_i$ and $T_j$), namely

\begin{equation*}    
P_{ij}^{\text{rad}}
\approx G_{ij}^{\text{rad}}(\bar T)\,(T_i-T_j),
\end{equation*}

where
\begin{equation*}
G_{ij}^{\text{rad}}(\bar T)
:= 4\,\sigma\,\bar T^{3}\,\epsilon_{ij}^\ast A_i F_{ij}.
\end{equation*}

\end{itemize}

These parameters cannot be fully determined analytically and require calibration based on measurements or high-fidelity simulations. Related data-driven approaches for modeling and correcting radiative thermal couplings in reduced-order spacecraft models have been explored \cite{benner2015survey, willcox2002balanced}.

In the following, radiative exchanges with the external environment that can be reliably estimated from design data are treated as known exogenous inputs. Conversely, unresolved or uncertain radiative interactions between components may be absorbed into the effective conductive coefficients to be identified. This choice is consistent with the interpretation of the model parameters as effective quantities rather than exact physical properties.

\subsection{Limitations of direct lumped reductions}

The parameters in Lumped Parameter Thermal Models (LPTMs) are often challenging to estimate from nominal material data and geometric specifications. Variations in material properties, uncertainties in geometry, and neglected interactions lead to differences between the idealized model assumptions and the observed effective behavior. The key assumption underlying our approach is that, despite these simplifications, the reduced-order model retains the fundamental structure of thermal physics: heat exchanges between lumped nodes are governed by temperature gradients, energy is conserved, and the dynamics remain consistent with thermodynamic principles. Under this assumption, the LPTM can be interpreted as a physics-constrained surrogate, where the governing equations define the functional form of the dynamics and the parameters are treated as effective lumped coefficients to be identified from observed thermal trajectories. Rather than deriving these coefficients from first principles, we learn them from data while enforcing the structural constraints imposed by the thermal network topology.
This motivates the formulation of parameter identification as a data-driven inverse problem, constrained by the physical structure of the thermal network.

\section{Inverse Parameter Identification Framework}
\label{sec:3}

The objective of the proposed framework is to identify a compact and lumped thermal model from partial observations, which naturally leads to a parameter estimation problem within the theory of inverse problems \cite{stuart2010inverse,kaipio2005inverse}. Rather than deriving directly the lumped model parameters as defined in Section~\ref{sec:lptm}, we reinterpret the lumped formulation introduced in Section 2 as a parameterized dynamical system whose coefficients are inferred directly from data.

\subsection{Effective Lumped Thermal Model}

Let $\phi := \{C_i, G_{ij}\}$ denote the physical parameters that would characterize an ideal lumped thermal model derived from the underlying heat-transfer problem. In practice these quantities are not directly accessible and cannot be reliably determined from geometry and material properties alone.
Instead, we introduce a set of effective coefficients $\theta = \{\gamma_i, \delta_{ij}\}$ that parameterize a reduced dynamical system reproducing the observed temperature evolution at the sensor level. These coefficients play the role of nodal (inverse) capacitances and edge conductances within the proxy model, but they should not be interpreted as exact physical quantities.

Let \(T(t) \in \mathbb{R}^n\) denote the vector of nodal temperatures.
The reduced thermal dynamics are described by the system

\begin{equation*}
\dot T(t) = f_\theta(T(t), u(t), t),
\end{equation*}

which can be written in the equivalent lumped form

\begin{equation*}
C(\theta)\dot T(t) = -L(\theta)T(t) + u(t).
\end{equation*}

Here \(C(\theta)\) collects the effective nodal capacitances and
\(L(\theta)\) is the weighted graph Laplacian associated with the
thermal network. At the node level, the dynamics can be written as

\begin{equation*}
\dot T_i =
\gamma_i
\left(
u_i +
\sum_{j \in N(i)} \delta_{ij}(T_j - T_i)
\right),
\end{equation*}

where \(N(i)\) denotes the set of neighbors of node \(i\). The coefficients \(\gamma_i\) and \(\delta_{ij}\) represent effective thermal properties that must be estimated from data. The coefficient \(\gamma_i\) can be interpreted as an effective inverse thermal capacitance, i.e., \(\gamma_i \approx C_i^{-1}\), while \(\delta_{ij}\) represents an effective thermal conductance between nodes \(i\) and \(j\).

\subsection{Supervised Learning Perspective}

The dynamical system introduced in Section 3.1 naturally inherits the graph structure of the underlying thermal network. Each temperature state corresponds to a node of the graph, while thermal couplings between components are represented by edges. Within this representation, the evolution of each node results from the combination of local power inputs and pairwise heat exchanges with neighboring nodes. The conductive interaction between two nodes $i$ and $j$ can be written as

\begin{equation*}
m_{ij}(T_i,T_j;\theta)
=
\delta_{ij}(T_j - T_i),
\end{equation*}

which represents the heat transfer along the edge connecting the two
nodes. The total contribution acting on node $i$ is obtained by
aggregating the interactions with all adjacent nodes,

\begin{equation*}
\sum_{j \in N(i)} m_{ij}(T_i,T_j;\theta).
\end{equation*}

The nodal dynamics therefore take the form

\begin{equation*}
\dot T_i =
\gamma_i
\left(
u_i +
\sum_{j \in N(i)} m_{ij}(T_i,T_j;\theta)
\right).
\end{equation*}

This formulation coincides with the message-passing structure commonly used in graph-based dynamical systems. The same functional form of the interaction law is applied across all edges of the network, while the parameters \(\delta_{ij}\) remain edge-specific and are learned from data. From this perspective, the lumped thermal model can be interpreted as a structured graph dynamical system in which the interaction functions are fixed by the physical heat transfer law, while the associated coefficients are unknown parameters to be identified.

This observation establishes a direct connection with the framework of Graph Neural Ordinary Differential Equations, which generalize Neural ODE models to structured systems \cite{chen2018neuralode,dupont2019augmented}. In that setting, the
state evolution is described by a parameterized vector field defined on a graph, and the unknown parameters are identified from trajectory data through gradient-based optimization.

Although the present model does not introduce neural network
parametrizations for the interaction functions, its structure matches that of a Graph Neural ODE with physically prescribed message functions. In contrast to fully data-driven approaches such as Physics-Informed Neural Networks, the proposed framework retains an explicit physical structure and focuses on parameter identification rather than function approximation \cite{raissi2019pinn,karniadakis2021piml}. This reinterpretation allows the calibration problem to be treated within the standard supervised learning setting commonly adopted for Neural ODE models.

\subsection{Inverse Problem Formulation}

Let $\{t_k\}_{k=0}^N$ denote the sampling times and 
$\{u_k, y_k\}_{k=0}^N$ the available dataset, where $u_k$
represents the known thermal inputs and $y_k$ the measured
temperatures at the sensor locations. Given a parameter vector $\theta$, the effective thermal model introduced in Section 3.1 defines a dynamical system whose evolution is governed by

\begin{equation*}
\dot T(t) = f_{\theta}(T(t),u(t),t).
\end{equation*}

For a prescribed initial condition $T(t_0)$ and input
trajectory $u(t)$, this system generates a temperature
trajectory $\{T_k(\theta)\}_{k=0}^{N}$, where we denote $T_k(\theta):=T(t_k;\theta)$. The identification task consists in determining the parameter vector $\theta$ that best reproduces the observed trajectories. This leads to a trajectory-matching optimization problem, consistent with classical system identification approaches based on time-series data \cite{ljung1994identifiability},

\begin{equation}
\theta^\star =
\arg\min_{\theta}
\sum_{k=0}^{N}
\|T_k(\theta) - y_k\|_2^2.
\label{eq:traj_loss}
\end{equation}

From the perspective introduced in Section 3.2, this problem can be interpreted as a supervised learning task on dynamical systems defined over a graph. The model parameters $\theta$ play the role of trainable weights, while the temperature trajectories represent the training targets. 

In practice, solving this optimization problem requires a number of additional design choices, including the parameterization of the coefficients, the strategy used to compute gradients of the trajectory loss, and the numerical structure adopted to enable efficient training. These aspects are addressed in the next section, where we describe the computational architecture used to make the inverse problem numerically tractable while preserving the physical structure of the thermal model.

The trajectory matching differs from standard regression settings, as the supervision is provided through the evolution of the system rather than pointwise targets. In practice, for computational and stability reasons, the trajectory matching loss will be evaluated over multiple shorter time windows, as described in Section 4. The proposed framework can therefore be interpreted as a structured neural differential model on graphs, where physical constraints such as positivity, symmetry, and dissipativity are enforced by construction, while retaining end-to-end trainability from data.
\section{Computational Architecture for Stable Calibration}
\label{sec:4}

The inverse problem formulated in Section 3 corresponds to the identification of the parameter vector $\theta$ of a graph dynamical system from observed trajectories. While conceptually simple, solving this optimization problem in practice requires several architectural choices to ensure numerical stability and efficient training. Stability considerations mirror known challenges in the training of Neural ODE models over long time horizons \cite{dupont2019augmented,kidger2022survey}.

In particular, trajectory-based optimization of Neural ODE systems can become unstable when parameters temporarily leave the physically meaningful regime or when gradients are computed over long time horizons. For this reason the calibration framework adopts a number of structural design choices that guarantee stability of both the dynamics and the optimization process.

\subsection{Parallel Trajectory Batching}
\label{sec:parallel}

The trajectory-based loss requires repeated integration of the dynamical system over multiple time intervals. For long telemetry sequences evaluating the loss over the entire trajectory at each optimization step would be computationally too expensive and may lead to unstable gradients, since prediction errors accumulate over long horizons.

To address this issue, the available time series is partitioned into multiple shorter temporal segments that are processed simultaneously during training. Let $N_B$ denote the number of batches and $T_b$ the temporal length of each segment. The dataset is partitioned as
\begin{equation*}
D_b = \{u_{b,k}, y_{b,k}\}_{k=0}^{N_b}
\end{equation*}

where each batch corresponds to a trajectory initialized from a
different time window of the telemetry sequence. Each trajectory segment corresponds to a contiguous time window extracted from the full dataset. In this work, non-overlapping segments are used, although overlapping windows could also be considered.

Instead of evolving each trajectory independently, the system is extended to a larger disconnected graph obtained by replicating the original thermal network $N_B$ times. If $L(\theta)$ denotes the graph Laplacian of the original network, the corresponding operator for the composite system takes a block-diagonal form and corresponding state vector is constructed by concatenating the temperature states of all batches,

\begin{equation*}
\bar L(\theta) =
\mathrm{diag}\left(L(\theta),\ldots,L(\theta)\right),
\qquad \bar T =
\left[T_1^\top, T_2^\top, \ldots, T_{N_B}^\top\right]^\top .
\end{equation*}

The resulting system evolves according to 
\begin{equation}
\dot{\bar T}=f_\theta(\bar T,\bar u,t),   
\end{equation}
where each subgraph evolves independently but shares $\theta$.

\smallskip

This formulation allows all trajectories to be propagated within a single forward pass of the graph-based Neural ODE solver. Since the graph operations involved in the vector field evaluation consist primarily of sparse matrix–vector products and local message computations along edges, the resulting workload can be executed efficiently using parallel computing architectures.

In addition, processing multiple trajectories simultaneously improves the conditioning of the optimization problem.
Gradient updates are aggregated across batches corresponding to different initial conditions and forcing scenarios, providing a more informative training signal and reducing the risk of overfitting to individual trajectory segments.

\subsection{Bounded Parameter Embedding}
\label{sec:bounded}

The physical parameters of the thermal network are collected in the vector $\theta = \{\gamma_i,\delta_{ij}\}$. These coefficients have two structural properties that make direct optimization very difficult. First, they are all strictly positive. Second, their admissible values may span several orders of magnitude, since they correspond to thermal processes occurring at different characteristic time scales.
Very small coefficients correspond to thermally inactive nodes or edges, effectively removing interactions from the network. Conversely, very large coefficients represent increasingly fast heat-transfer processes. The high regime is problematic both numerically and conceptually: extremely large parameters introduce stiff dynamics that require very small integration timesteps, while from a modeling perspective they correspond to resolving thermal processes that fall outside the intended resolution of the lumped approximation.

To enforce physical constraints such as positivity and boundedness of the parameters, we introduce a smooth mapping from unconstrained variables \(\tilde{p} \in \mathbb{R}\) to physically admissible parameters \(p\).

Specifically, we define
\[
p = \Psi(\tilde{p}) = \frac{s_p}{2} \left( \tanh(2\tilde{p}) + 1 \right),
\]
where \(s_p > 0\) is a characteristic scale representing the maximum admissible value of the parameter.

This transformation ensures that \(p \in (0, s_p)\) for all \(\tilde{p} \in \mathbb{R}\), while preserving smoothness and enabling gradient-based optimization in the unconstrained space. The choice of the hyperbolic tangent provides symmetric saturation and avoids excessively small gradients near the origin, thereby improving numerical stability during training. Similar constrained parameterizations are commonly adopted to enforce positivity and boundedness in neural differential models \cite{drgona2020constrained,sanchez2020learning}. The scale \(s_p\) is selected based on prior physical knowledge or expected orders of magnitude of the parameters.

\subsection{Bounded Evolution During Training}
\label{sec:evo}

Even with the bounded parameterization introduced in Section~\ref{sec:bounded}, the optimization process may temporarily explore parameter values far from those corresponding to a stable thermal network. In such situations the predicted trajectories can diverge rapidly during the numerical integration of the dynamical system. This issue is particularly critical when the trajectory loss is evaluated through long forward rollouts, since unstable intermediate dynamics may lead to exploding states and corrupted gradient estimates.

To address this issue, we introduce a stabilization strategy that limits the range of admissible temperatures during the forward propagation of the model. Specifically, after each integration step of the unrolled dynamics, a smooth limiting function is applied to the temperature states. The transformation is designed so that it does not modify the dynamics within the physically relevant operating range, while providing bounded dynamics when the trajectory moves outside this region.

7uLet $T$ denote a generic nodal temperature. The limited value $\hat T$ is defined through a piecewise smooth transformation constructed using hyperbolic tangent functions

\begin{equation}
\begin{cases}
    v\tanh\left(T/v-1\right)+1v, &T < v\\
    T, &v \le T \le 2v\\
    v\tanh\left(T/v-2\right)+2v, &2v < T
\end{cases}
\label{eq:smoothing}
\end{equation}

The mapping is linear within the nominal operating interval $200 \,\mathrm{K} \le T \le 400 \,\mathrm{K}$ using $v=200\,\mathrm{K}$, so that temperatures in this range remain unchanged. Outside this interval the function transitions smoothly toward asymptotic bounds using hyperbolic tangent branches. In particular, the transformation approaches the lower asymptote at $0\,\mathrm{K}$ and the upper
asymptote at $600\,\mathrm{K}$.

This choice ensures that the transformation behaves as the identity in the physically meaningful temperature range, while preventing unbounded growth of the state variables during unstable intermediate phases of the optimization.

The temperature limiter acts only outside the physically relevant operating range and is primarily used to prevent numerical instabilities during the early stages of training, consistent with the numerical stability issues discussed in \cite{massaroli2021dissecting}, when parameters may take unrealistic values. As training progresses, the learned dynamics remain within the linear region of the limiter, making its effect negligible. In all reported experiments, the limiter is removed during evaluation, ensuring that the final model corresponds to the original unconstrained dynamics.

\subsection{Energy Conservation and Dissipation}

Beyond the numerical stabilization strategies introduced in the previous sections, the thermal model must preserve the structural properties of passive heat-transfer networks. In particular, conductive couplings between nodes are symmetric and cannot generate energy internally. Symmetry of conductive exchanges is enforced by construction through $\delta_{ij} = \delta_{ji}$, while nodal coefficients satisfy $\gamma_i>0$. Under these conditions the associated graph Laplacian $L(\theta)$ is symmetric positive semidefinite.

\smallskip 

Two energy-related quantities can be associated with the thermal network. First, the total thermal content of the system is defined as

\begin{equation}
H(T) = \sum_i C_i T_i = \mathbf{1}^\top C T .
\end{equation}

This quantity represents the total thermal content of the system relative to a chosen reference temperature. In the absence of external inputs, conductive exchanges only redistribute heat between nodes and therefore preserve this quantity. Indeed, using

\begin{equation}
C(\theta)\dot T = -L(\theta)T,
\end{equation}

one obtains

\begin{equation}
\dot H = \mathbf{1}^\top C(\theta)\dot T
       = -\mathbf{1}^\top L(\theta)T = 0,
\end{equation}

since the Laplacian satisfies

\begin{equation}
L(\theta)\mathbf{1}=0.
\end{equation}

In addition to the conserved thermal content, it is convenient to introduce the quadratic functional

\begin{equation}
E(T) = \frac12 T^\top C(\theta)T .
\end{equation}

Although not representing the physical thermal energy, this quantity acts as a Lyapunov function for the conductive dynamics.  Differentiating along system trajectories yields

\begin{equation}
\dot E
= T^\top C(\theta)\dot T
= -T^\top L(\theta)T .
\end{equation}

Using the structure of the graph Laplacian, the quadratic form can be written explicitly as

\begin{equation}
T^\top L(\theta)T =
\frac12 \sum_{i,j} \delta_{ij}(T_i-T_j)^2 ,
\end{equation}

which is always non-negative since $\delta_{ij} > 0$. Such properties are consistent with classical passivity-based formulations of networked physical systems.

Consequently, $\dot E \le 0$. Equality holds only when all connected nodes have identical temperatures, corresponding to a state of internal thermal equilibrium with no conductive heat fluxes.

Embedding symmetry and positivity directly in the parameterization therefore guarantees that the calibrated model preserves the passive and dissipative nature of thermal conduction throughout the entire training process. This ensures that the identified model remains within the class of passive and dissipative systems. Structure-preserving formulations of dynamical systems similarly aim to enforce consistency with underlying physical principles in learned models, including energy conservation and dissipation properties \cite{greydanus2019hnn,zhong2020symplectic}.
\section{Numerical Results}
\label{sec:5}

\begin{figure}[H]
\centering
\includegraphics{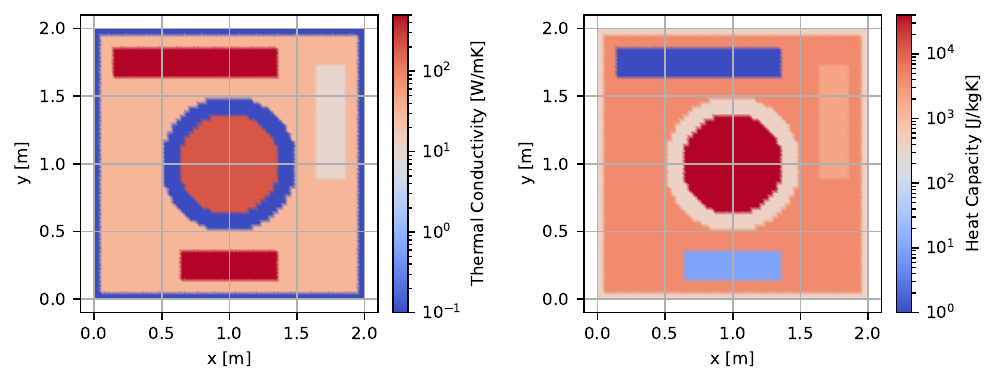}
\caption{Computational domain $\Omega = [0,2]^2$ m$^2$. Thermal conductivity field (left)
and volumetric heat capacity field (right).}\label{fig:domain}
\end{figure}

We first describe the test cases used to assess the proposed hybrid calibration framework. The computational domain represents a simplified yet physically meaningful cross-section of a satellite thermal subsystem. Let $\Omega \subset \mathbb{R}^2$ denote the spatial domain under consideration. Specifically, we consider a square two-dimensional domain $\Omega$, whose boundary $\partial \Omega$ is surrounded by a strongly insulating layer, mimicking the presence of multi-layer insulation and minimizing heat exchange with the external environment. Within $\Omega$, four disjoint regions $\{\Omega_i\}_{i=1}^4$ corresponding to internal subsystems are considered; three of these subsystems, $\Omega_1,\Omega_2,\Omega_3$, are monitored and have different thermal properties, while $\Omega_4$ is thermally isolated and unmonitored. Each region $\Omega_i$ is assumed homogeneous, with constant thermal conductivity $\lambda_i$ and volumetric heat capacity $\rho_i c_i$, as shown in  Figure~\ref{fig:domain}

Ground truth data are generated by solving the transient heat conduction problem for the temperature field $T(x,t)$ on $\Omega$ using a Finite Element Method (FEM), a standard discretization approach for heat transfer problems. 
The temperature field $T(x,t)$ is governed by:
\begin{equation}
c_f \, \partial_t T - \nabla \cdot (\lambda_f \nabla T)
= \sum_j f_j(x,t)\, \chi_{\Omega_j}(x)
\quad \text{in } \Omega,
\end{equation}
where $c_f$ and $k_f$ denote the thermal capacity and conductivity, respectively. The terms $f_j(x,t)$ represent localized volumetric heat sources acting on subdomains $\Omega_j \subset \Omega$. On the boundary $\partial\Omega = \cup_{i=1}^4 \Gamma_i$, nonlinear Neumann boundary conditions are imposed:
\begin{equation}
- \lambda_f \nabla T \cdot n
= q_i(t) + \sigma \varepsilon z_{\mathrm{len}}
\left(T ^4 - T_{\mathrm{ext}}^4 \right)
\quad \text{on } \Gamma_i ,
\end{equation}
accounting for surface heat fluxes $q_i(t)$ and radiative heat losses linearized towards an external temperature $T_{\mathrm{ext}}$.

\subsection{Weak Formulation and Data Generation}
\label{sec:data}

The synthetic dataset is generated by solving the two-dimensional heat equation with volumetric sources, prescribed boundary heat fluxes, and radiative cooling. Using an IMEX scheme for time discretization with timestep $\Delta t$, the weak formulation reads as follows: find $u \in H^1(\Omega)$ such that for all test functions $v \in H^1(\Omega)$

\begin{equation*}
\int_\Omega c_f T v \, dx 
+ \Delta t \int_\Omega \lambda_f \nabla T \cdot \nabla v \, dx 
- \Delta t \sum_{i=1}^{4} \int_{\Gamma_i} 
\left( \sigma \varepsilon z_{\text{len}} T_{\text{old}}^{3} \right) T v \, ds
=
\int_\Omega c_f T_{\text{old}} v \, dx 
+ \Delta t \sum_{i=1}^{4} q_i \int_{\Gamma_i} v \, ds
+ \Delta t \sum_{j} f_j \int_{\Omega_j} v \, dx.
\end{equation*}

Here $T$ denotes the temperature at the current time step and $T_{\text{old}}$ the temperature at the previous one. The terms $q_i$ correspond to prescribed heat fluxes on the boundary segments $\Gamma_i$, while $f_j$ denote volumetric heat sources acting on subdomains $\Omega_j \subset \Omega$. Radiative cooling is modeled through a boundary contribution derived from the Stefan--Boltzmann law. The nonlinear $T^4$ dependence is linearized around the previous state and assuming $T_{ext}\approx0K$, leading to the coefficient $\sigma \varepsilon z_{\text{len}} T_{\text{old}}^{3}$ and yielding a stable fully implicit formulation.

The computational domain is $\Omega = [0,2] \times [0,2]$ m$^2$ and is discretized with a structured $64 \times 64$ mesh of rectangular elements. Continuous Galerkin finite elements of degree one (P1) are employed, resulting in $65 \times 65 = 4225$ degrees of freedom. Time integration is  solved by direct LU factorization.

\subsection{Generation of Time Series and Heat Sources}

The resulting FEM solution serves as a surrogate for ground-truth physics, from which synthetic measurements are extracted. All results presented in this work are obtained using simulated data in a controlled setting \ref{sec:data}, where the thermal dynamics, material properties, and boundary conditions are known. High-fidelity FEM simulations provide a reference solution that is not available in real telemetry, allowing us to quantify to what extent the proposed framework is able to recover a consistent lumped thermal model from sparse temperature measurements.

To emulate a realistic operational scenario with sparse sensing, only a limited number of temperature sensors is assumed to be available. All sensors provide pointwise temperature measurements and constitute the only observable data used by the proposed model. From each sensor location $\{\bm{x}_k\}_{k=0}^8$, a discrete temperature time series
$\{T(\bm{x}_k, t_n)\}_{n=0}^{N_t}$ is extracted, and the collection of these sparse
observations forms the input dataset for the hybrid model.

In these test cases, 4 sensors are located on the boundary of the square domain, one on each side, and are referred to as boundary nodes. These sensors define the interface through which radiative heat exchange between the corresponding satellite face and deep space is applied. Accordingly, the boundary nodes are endowed with the known geometric and optical properties of the associated satellite face, such as radiating area and emissivity. The remaining sensors are placed inside the monitored subdomains, with some subdomains containing more than one sensor.

\begin{figure}[H]
\centering
\includegraphics{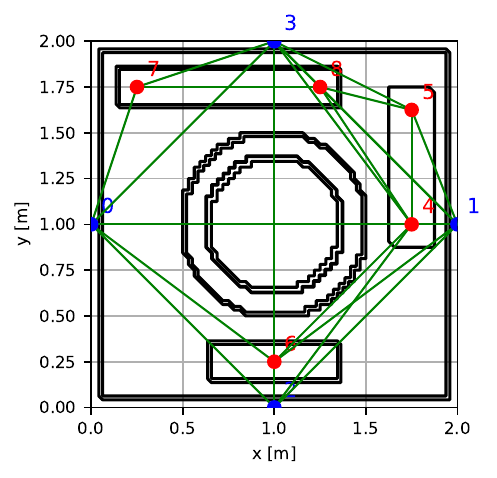}
\caption{Thermal graph on the
computational domain. Boundary nodes (1--4) handle radiative exchange with deep space;
internal nodes (5--8) correspond to sensor locations within the monitored subdomains.
Green edges denote conductive couplings.}\label{fig:schematics}
\end{figure}

The thermal forcing time series were generated using two methodologies to represent deterministic and stochastic phenomena. Periodic forcings, simulating solar irradiation, were modeled through sinusoidal functions with characteristic orbital period, allowing to reproduce the solar exposure and eclipse phases typical of low Earth orbit missions. For stochastic forcings, representing both heating lines and thermal dissipations from internal electronic devices, a Gaussian process approach with radial basis function kernel was employed. This choice is motivated by the intrinsically variable and unpredictable nature of internal thermal loads, which depend on the non-deterministic activation of electronic components, payload operational cycles, and variable usage conditions during the mission.  In the test below, Tests A, B, and C refer to configurations with different characteristic periods of the internal sources, respectively a correlation length of 700, 840 and 560.

\subsection{Training Strategy}
\label{sec:training}

The training and validation dataset is processed in parallel according to the block-diagonal formulation of
Section~\ref{sec:parallel}.
Batch length $B$ is chosen to be sub-orbital, i.e.\ shorter than the characteristic thermal period of
approximately $100$ minutes.
The sensitivity of the results to the choice of $T_b$ is later investigated in
Section~\ref{sec:resultsB}.

\smallskip

To prevent numerical divergence during early optimization phases, we apply a smooth saturation \eqref{eq:smoothing} to nodal temperatures at each integration step. 
The transformation coincides with the identity within the physically admissible operating range $[200\,\mathrm{K},\,400\,\mathrm{K}]$, and smoothly compresses values outside this interval through a differentiable mapping with bounded asymptotes. 
As training progresses and the identified parameters enter a physically meaningful regime, predicted temperatures remain within this interval and the saturation acts as the identity. 
Consequently, it does not alter the learned dynamics in the operational regime and serves as a stabilization mechanism during early training.

\smallskip

Model parameters are optimized by minimizing the trajectory-based loss defined in ~\eqref{eq:traj_loss}. Each run is trained for a maximum of $5000$ epochs, with early stopping applied when the validation loss ceases to decrease over a fixed patience window. Parameters are optimized using an adaptive first-order gradient method (ADAM \cite{kingma2015adam}), applied to the ODE-constrained trajectory mismatch functional, without learning rate scheduler; runs are conducted independently at fixed learning rates of $10^{-2}$ and $10^{-3}$ to explore the sensitivity of convergence to this hyperparameter.

\smallskip

Model performance is evaluated on the long-horizon test segment $\mathcal{D}_{\mathrm{test}}$, which spans several complete orbital cycles and two metrics are reported: the root mean squared error~(RMSE) between predicted and reference temperature trajectories and the Pearson correlation coefficient~(PCC) \eqref{eq:pcc}, which quantifies the linear correlation between predicted and observed time series independently of amplitude offsets.

\begin{equation}
PCC = \frac{\text{cov}(X, Y)}{\sigma_X \cdot \sigma_Y} = \frac{\sum_{i=1}^{n}(x_i - \bar{x})(y_i - \bar{y})}
         {\sqrt{\sum_{i=1}^{n}(x_i - \bar{x})^2 \cdot \sum_{i=1}^{n}(y_i - \bar{y})^2}}
\label{eq:pcc}
\end{equation}

\subsection{Computational complexity}
Let $G = (V,E)$ denote the thermal graph. 
Each evaluation of the vector field $f_\theta$ consists of local message computations along edges followed by node-wise aggregation. 
Therefore, the computational cost per timestep scales as
\[
\mathcal{O}(|E| + |V|),
\]
since each edge contributes a constant-cost interaction and each node performs a constant-cost update. 

Under the block-diagonal batching strategy of Section~ \ref{sec:parallel}, $B$ independent trajectories are propagated in parallel, leading to an overall cost 
\[
\mathcal{O}\big(B(|E| + |V|)\big)
\]
per timestep, while sharing the same parameter vector $\phi$. 
Memory requirements scale linearly with $B(|E| + |V|)$.

\subsection{Test Case Configurations}
\label{sec:resultsA}

Three test configurations are considered to assess the proposed framework under progressively complex thermal scenarios. The first (Table~\ref{tab:resultsE}) involves only external radiative forcing applied at the domain boundary, representing orbital heat fluxes in the absence of internal dissipation.
The second (Table~\ref{tab:resultsS}) features two internal volumetric heat sources acting on distinct subdomains, mimicking the activation of heater lines used in practice to maintain onboard components within their safe operating temperature ranges.
The third configuration (Table~\ref{tab:resultsM}) combines both external and internal forcings simultaneously and constitutes the reference benchmark for all subsequent analyses reported in Section~\ref{sec:resultsB}.

\smallskip

Convergence is consistently achieved within a small number of iterations, often well below the maximum allowed. This is due to to three structural properties of the proposed framework. First, the strong physical constraints embedded in the neural field architecture (positivity, symmetry, and energy conservation) substantially reduce the effective dimensionality of the optimization landscape by confining the search to the physically admissible parameter space from the first iteration.
Second, the bounded parameterization introduced in Section~\eqref{sec:bounded} prevents gradient explosion and ensures smooth differentiability throughout training. Third, the parallelization strategy described in Section~\ref{sec:parallel} provides gradient updates averaged over multiple initial conditions and input scenarios simultaneously, yielding a more informative and stable training signal per iteration compared to single-trajectory
optimization.

\smallskip

\begin{figure}[H]
    \centering
    \includegraphics[clip]{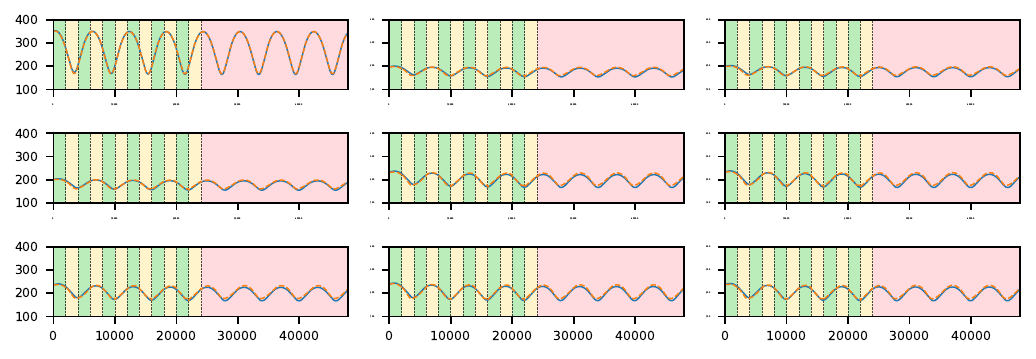}
    \caption{Predicted (orange) and reference (blue) temperature trajectories for the external
source configuration. }

    \label{fig:result1}
\end{figure}

\begin{figure}[H]
    \centering
    \includegraphics[clip]{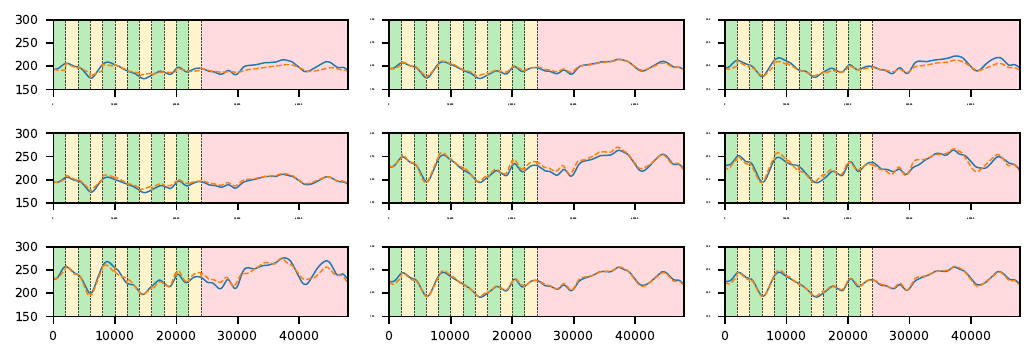}
    \caption{Predicted (orange) and reference (blue) temperature trajectories for the internal heater line configuration.}
    \label{fig:result2}
\end{figure}

\begin{figure}[H]
    \centering
    \includegraphics[clip]{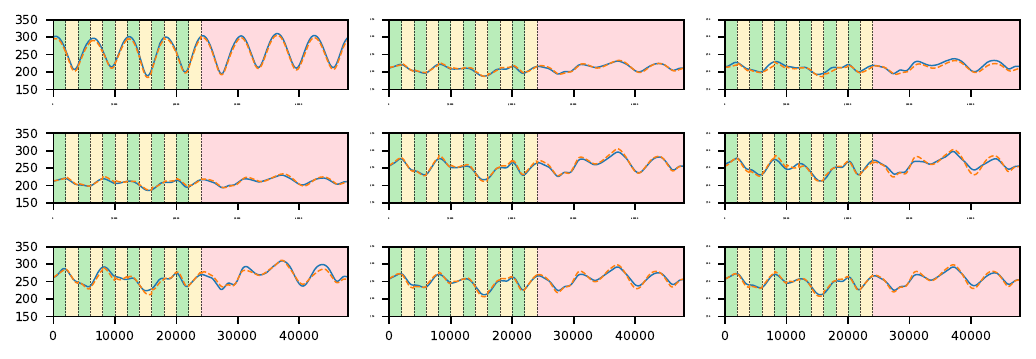}
    \caption{Predicted (orange) and reference (blue) temperature trajectories for the benchmark test, combining external and internal forcings.}
    \label{fig:result3}
\end{figure}

Representative results for the three configurations are reported in Figures~\ref{fig:result1}-\ref{fig:result2}-\ref{fig:result3}. Each panel displays the predicted and reference temperature trajectories over the full time horizon considered. The background shading identifies the three data partitions: the green region corresponds to the training segments, the yellow region to the validation segments, and the red region to the long-horizon test set $\mathcal{D}_{\mathrm{test}}$. In each panel, the reference FEM temperature is shown in blue and the model prediction in orange. 

\begin{table}[H]
\centering
\begin{tabular}{lllrrrrr}
\toprule
Dataset & Iterations & Max Iter & Learning Rate & Train Loss & Valid Loss & Test RMSE & Test PCC \\
\midrule
Test A & 2000 & 2000 & $1.00 \times 10^{-2}$ & $1.43 \times 10^{-2}$ & $1.66 \times 10^{-2}$ & $2.66 \times 10^{-2}$ & $9.98 \times 10^{-1}$ \\
Test A & 2073 & 5000 & $1.00 \times 10^{-2}$ & $1.43 \times 10^{-2}$ & $1.65 \times 10^{-2}$ & $2.66 \times 10^{-2}$ & $9.98 \times 10^{-1}$ \\
Test A & 5000 & 5000 & $1.00 \times 10^{-3}$ & $1.50 \times 10^{-2}$ & $1.75 \times 10^{-2}$ & $2.71 \times 10^{-2}$ & $9.98 \times 10^{-1}$ \\
\bottomrule
\end{tabular}
\caption{Calibration results for the external forcing configuration.}
\label{tab:resultsE}
\end{table}

\begin{table}[H]
\centering
\begin{tabular}{lllrrrrr}
\toprule
Dataset & Iterations & Max Iter & Learning Rate & Train Loss & Valid Loss & Test RMSE & Test PCC \\
\midrule
Test A & 2000 & 2000 & $1.00 \times 10^{-2}$ & $2.00 \times 10^{-2}$ & $2.03 \times 10^{-2}$ & $2.16 \times 10^{-2}$ & $9.79 \times 10^{-1}$ \\
Test A & 4210 & 5000 & $1.00 \times 10^{-2}$ & $1.93 \times 10^{-2}$ & $1.96 \times 10^{-2}$ & $2.08 \times 10^{-2}$ & $9.81 \times 10^{-1}$ \\
Test A & 5000 & 5000 & $1.00 \times 10^{-3}$ & $2.31 \times 10^{-2}$ & $2.31 \times 10^{-2}$ & $2.41 \times 10^{-2}$ & $9.75 \times 10^{-1}$ \\
Test B & 2000 & 2000 & $1.00 \times 10^{-2}$ & $2.08 \times 10^{-2}$ & $1.99 \times 10^{-2}$ & $2.40 \times 10^{-2}$ & $9.73 \times 10^{-1}$ \\
Test B & 4195 & 5000 & $1.00 \times 10^{-2}$ & $1.96 \times 10^{-2}$ & $1.88 \times 10^{-2}$ & $2.15 \times 10^{-2}$ & $9.79 \times 10^{-1}$ \\
Test B & 5000 & 5000 & $1.00 \times 10^{-3}$ & $2.27 \times 10^{-2}$ & $2.25 \times 10^{-2}$ & $2.79 \times 10^{-2}$ & $9.64 \times 10^{-1}$ \\
Test C & 2000 & 2000 & $1.00 \times 10^{-2}$ & $1.60 \times 10^{-2}$ & $2.25 \times 10^{-2}$ & $3.96 \times 10^{-2}$ & $9.38 \times 10^{-1}$ \\
Test C & 4753 & 5000 & $1.00 \times 10^{-2}$ & $1.40 \times 10^{-2}$ & $2.30 \times 10^{-2}$ & $3.44 \times 10^{-2}$ & $9.50 \times 10^{-1}$ \\
Test C & 5000 & 5000 & $1.00 \times 10^{-3}$ & $1.82 \times 10^{-2}$ & $3.33 \times 10^{-2}$ & $3.98 \times 10^{-2}$ & $9.34 \times 10^{-1}$ \\
\bottomrule
\end{tabular}
\caption{Calibration results for the internal heater line configuration.}
\label{tab:resultsS}
\end{table}

\begin{table}[H]
\centering
\begin{tabular}{lllrrrrr}
\toprule
Dataset & Iterations & Max Iter & Learning Rate & Train Loss & Valid Loss & Test RMSE & Test PCC \\
\midrule
Test A & 2000 & 2000 & $1.00 \times 10^{-2}$ & $1.74 \times 10^{-2}$ & $2.18 \times 10^{-2}$ & $2.03 \times 10^{-2}$ & $9.86 \times 10^{-1}$ \\
Test A & 2456 & 5000 & $1.00 \times 10^{-2}$ & $1.72 \times 10^{-2}$ & $2.12 \times 10^{-2}$ & $1.97 \times 10^{-2}$ & $9.87 \times 10^{-1}$ \\
Test A & 5000 & 5000 & $1.00 \times 10^{-3}$ & $1.93 \times 10^{-2}$ & $2.43 \times 10^{-2}$ & $2.23 \times 10^{-2}$ & $9.82 \times 10^{-1}$ \\
Test B & 2000 & 2000 & $1.00 \times 10^{-2}$ & $1.77 \times 10^{-2}$ & $1.91 \times 10^{-2}$ & $2.13 \times 10^{-2}$ & $9.85 \times 10^{-1}$ \\
Test B & 2336 & 5000 & $1.00 \times 10^{-2}$ & $1.76 \times 10^{-2}$ & $1.90 \times 10^{-2}$ & $2.11 \times 10^{-2}$ & $9.85 \times 10^{-1}$ \\
Test B & 5000 & 5000 & $1.00 \times 10^{-3}$ & $1.92 \times 10^{-2}$ & $2.07 \times 10^{-2}$ & $2.37 \times 10^{-2}$ & $9.82 \times 10^{-1}$ \\
Test C & 2000 & 2000 & $1.00 \times 10^{-2}$ & $1.53 \times 10^{-2}$ & $2.18 \times 10^{-2}$ & $2.09 \times 10^{-2}$ & $9.85 \times 10^{-1}$ \\
Test C & 5000 & 5000 & $1.00 \times 10^{-2}$ & $1.28 \times 10^{-2}$ & $1.97 \times 10^{-2}$ & $1.93 \times 10^{-2}$ & $9.88 \times 10^{-1}$ \\
Test C & 5000 & 5000 & $1.00 \times 10^{-3}$ & $1.73 \times 10^{-2}$ & $2.43 \times 10^{-2}$ & $2.32 \times 10^{-2}$ & $9.81 \times 10^{-1}$ \\
\bottomrule
\end{tabular}
\caption{Calibration results for the benchmark configuration.}
\label{tab:resultsM}
\end{table}

\begin{figure}[htbp]
\centering

\begin{subfigure}[t]{0.30\textwidth}
    \centering
    \includegraphics{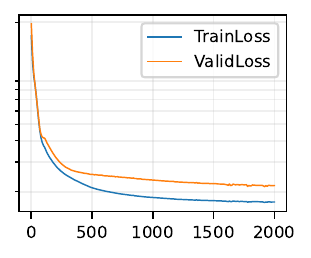}
    \caption{Training and validation loss versus epoch}
    \label{fig:sub1}
\end{subfigure}
\hfill
\begin{subfigure}[t]{0.30\textwidth}
    \centering
    \includegraphics{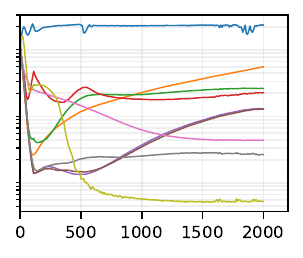}
    \caption{Evolution of the learned nodal capacitances $\gamma_i$}
    \label{fig:sub2}
\end{subfigure}
\hfill
\begin{subfigure}[t]{0.30\textwidth}
    \centering
    \includegraphics{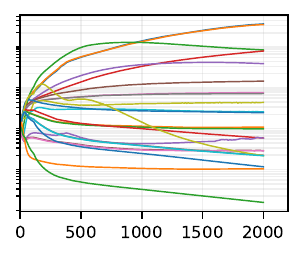}
    \caption{Evolution of the learned edge conductances $\delta_{ij}$.}
    \label{fig:sub3}
\end{subfigure}

\caption{Training diagnostics for the benchmark configuration.}
\label{fig:evolution}
\end{figure}

The convergence history and the evolution of the learned parameters over training iterations are reported in Figure~\ref{fig:evolution} for the benchmark configuration combining external radiative forcing and internal heater line activations. The training and validation losses decrease monotonically and reach a stable plateau within a limited number of iterations, consistent with the early stopping counts reported in Table~\ref{tab:resultsM}. The learned physical parameters converge smoothly to well-defined values without oscillations or instabilities, confirming that the bounded parameterization of Section~\ref{sec:bounded} effectively conditions the optimization landscape. An analogous convergence behavior is observed across all other test configurations, indicating that the framework is robust with respect to the thermal forcing.

\subsection{Sensitivity Analysis and Robustness Assessment}
\label{sec:resultsB}

All analyses reported in this section are conducted on the benchmark configuration introduced in Section~\ref{sec:resultsA}, which combines simultaneous external radiative forcing and internal heater line activations. Starting from this common reference setting, the training configuration and the properties of the input data are systematically varied in order to assess how each factor affects the behavior of the calibrated model.

\subsubsection*{Effect of Training Window Length}

In this experiment, the total amount of data available for training and validation is kept fixed, while the number of batches $B$ is varied. Since the total dataset length is held constant, increasing $B$ directly reduces the temporal length $T_b$ of each individual training window, and vice versa. This setup allows the effect of the window length on the identification process to be isolated without altering the overall information content of the dataset. The numerical results are summarized in Table~\ref{tab:test-window}.

When using the trajectory-based loss, each window is integrated from a single initial condition, so that prediction errors accumulate over the rollout horizon. Longer windows therefore introduce a more challenging optimization problem in the early training stages, reflected in slower convergence and in a computational cost per epoch that grows approximately linearly with $T_b$, as confirmed by the training times reported in Table~\ref{tab:test-window}. At the opposite extreme, very short windows provide only local supervision, leaving the model insufficiently constrained by long-range thermal dynamics and exposing it to overfitting, as evidenced by the markedly higher test error obtained at $T_b = 500$.

Beyond a minimum threshold, however, further increases in window length yield diminishing returns in predictive accuracy while substantially increasing computational cost. Results show that an intermediate window length provides the best trade-off between long-horizon generalization and computational efficiency. The optimal choice can therefore be identified as the shortest window that still guarantees stable performance on the test set.

\begin{table}[H]
\centering
\begin{tabular}{rccccr}
\toprule
$T_b$ & Train Loss & Valid Loss & Test RMSE & Test PCC & Training Time \\
\midrule
 500 & $1.67 \times 10^{-2}$ & $1.63 \times 10^{-2}$ & $5.05 \times 10^{-2}$ & $9.04 \times 10^{-1}$ & 5m\,56s \\
1000 & $1.72 \times 10^{-2}$ & $1.82 \times 10^{-2}$ & $2.00 \times 10^{-2}$ & $9.87 \times 10^{-1}$ & 10m\,26s \\
2000 & $1.74 \times 10^{-2}$ & $2.18 \times 10^{-2}$ & $2.03 \times 10^{-2}$ & $9.86 \times 10^{-1}$ & 21m\,30s \\
4000 & $1.97 \times 10^{-2}$ & $2.39 \times 10^{-2}$ & $1.92 \times 10^{-2}$ & $9.87 \times 10^{-1}$ & 35m\,54s \\
\bottomrule
\end{tabular}
\caption{Effect of the number of training batches $B$ on model performance and computational
cost. Total dataset length is held constant.}
\label{tab:test-window}
\end{table}

\begin{figure}[htbp]
\centering

\begin{subfigure}[t]{0.48\textwidth}
    \centering
    \includegraphics[clip]{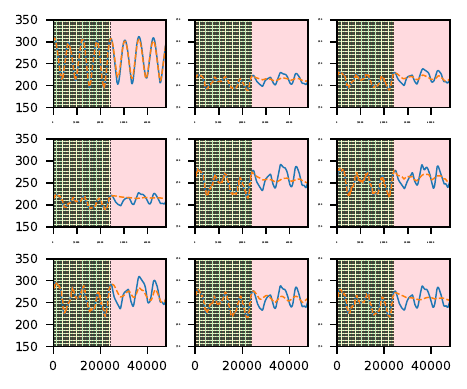}
    \caption{Result $T_b=500$}
\end{subfigure}
\hfill
\begin{subfigure}[t]{0.48\textwidth}
    \centering
    \includegraphics[clip]{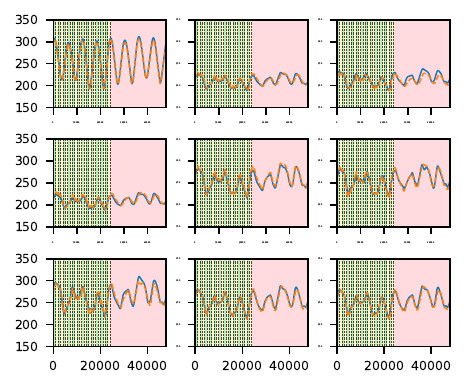}
    \caption{Result $T_b=1000$}
\end{subfigure}

\begin{subfigure}[t]{0.48\textwidth}
    \centering
    \includegraphics[clip]{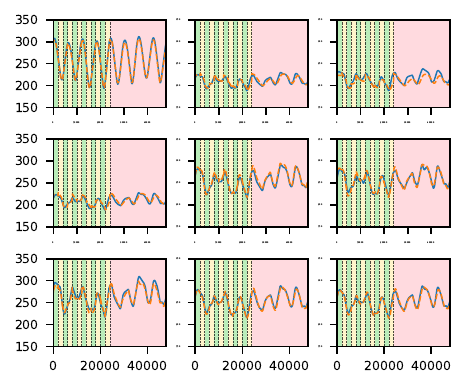}
    \caption{Result $T_b=2000$}
\end{subfigure}
\hfill
\begin{subfigure}[t]{0.48\textwidth}
    \centering
    \includegraphics[clip]{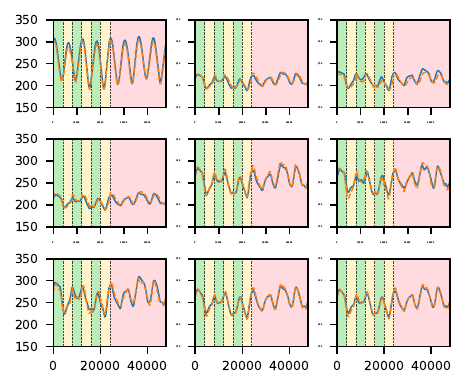}
    \caption{Result $T_b=4000$}
\end{subfigure}

\caption{Effect of training window length on the predicted temperature trajectories. Each panel corresponds to a different value of $T_b$, with background shading identifying training (green), validation (yellow), and test (red) segments.}
\label{fig:window}
\end{figure}

\newpage

\subsubsection*{Noise Robustness}
\label{sec:noise}

In this experiment, the robustness of the identification procedure is assessed in the presence of measurement noise. Gaussian white noise is added both to the initial conditions and to the temperature labels used for training and validation. Two noise levels are considered, equal to $1\%$ and $5\%$ of the mean temperature value $\overline{x}$. More precisely, at each time step a perturbation is sampled from a normal distribution with zero mean and standard deviation equal to $0.01\,\overline{x}$ and $0.05\,\overline{x}$, respectively. The corresponding results are reported in Table~\ref{tab:noise_results} alongside the noiseless baseline.

In the $1\%$ noise regime, model performance remains essentially unchanged with respect to the noiseless case. The trajectory loss, test error, and PCC values are all comparable to the baseline, indicating that the calibration framework is not sensitive to small-amplitude perturbations. When the noise level is increased to $5\%$, the test error increases accordingly but remains bounded and consistent with the noise magnitude, suggesting that the residual discrepancy is dominated by label uncertainty rather than model mismatch. Importantly, no divergence or instability is observed during training or autoregressive rollout in either regime.

To assess whether errors accumulate over the test horizon, the test segment is partitioned into sequential windows and the distribution of windowed errors is analyzed through box plots. The resulting  distributions remain statistically stable across all segments in both noise scenarios, confirming that the learned dynamics do not amplify measurement noise and that long-horizon robustness is preserved even under perturbed observations.

\begin{figure}[H]
\centering

\begin{subfigure}[t]{0.48\textwidth}
    \centering
    \includegraphics{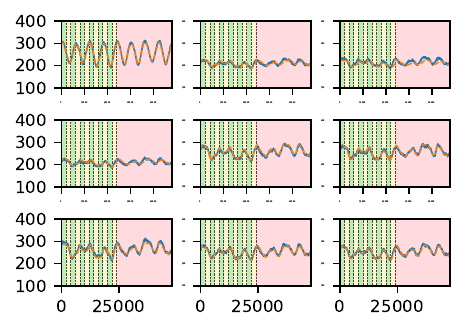}
    \caption{Result with $1\%$ Noise}
\end{subfigure}
\hfill
\begin{subfigure}[t]{0.48\textwidth}
    \centering
    \includegraphics{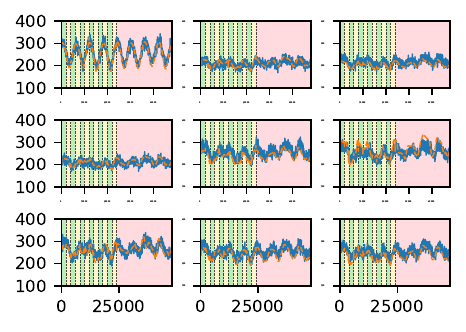}
    \caption{Result with $5\%$ Noise}
\end{subfigure}

\vspace{0.5em}

\begin{subfigure}[t]{0.48\textwidth}
    \centering
    \includegraphics{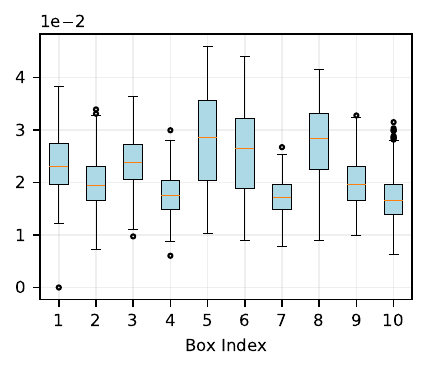}
    \caption{Box-Plot with $1\%$ Noise}
\end{subfigure}
\hfill
\begin{subfigure}[t]{0.48\textwidth}
    \centering
    \includegraphics{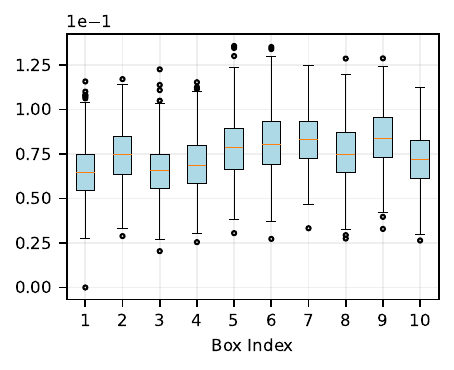}
    \caption{Box-Plot with $5\%$ Noise}
\end{subfigure}

\caption{Noise robustness results. Top: predicted (orange) vs.\ reference (blue) trajectories
under $1\%$ (left) and $5\%$ (right) measurement noise. Bottom: box plots of windowed test
errors for each noise level, confirming the absence of error accumulation over long horizons.}
\label{fig:noise}
\end{figure}

\begin{table}[ht]
\centering
\begin{tabular}{lccccc}
\hline
Noise & Train Error & Valid Error & Test RMSE & Test PCC & Time Total \\
\hline
$0\%$  & $1.740 \times 10^{-2}$ & $2.177 \times 10^{-2}$ & $2.028 \times 10^{-2}$ & $9.855 \times 10^{-1}$ & 21m\,30s \\
$1\%$  & $2.013 \times 10^{-2}$ & $2.319 \times 10^{-2}$ & $2.201 \times 10^{-2}$ & $9.832 \times 10^{-1}$ & 18m\,15s \\
$5\%$  & $5.757 \times 10^{-2}$ & $6.045 \times 10^{-2}$ & $5.679 \times 10^{-2}$ & $8.960 \times 10^{-1}$ & 18m\,43s \\
\hline
\end{tabular}
\caption{Performance metrics under different levels of measurement noise added to training labels
and initial conditions.}
\label{tab:noise_results}
\end{table}

\subsection{Training Stability Analysis}

A reliable calibration framework must produce consistent results across repeated training runs, even in the presence of measurement noise. In this section, training stability is evaluated from two complementary perspectives: the reproducibility of the performance metrics and the consistency of the learned physical parameters across runs.

The benchmark configuration is trained ten times on independently generated datasets, each corrupted by Gaussian white noise at a fixed level of $1\%$ of the mean temperature value $\overline{x}$ as described in Section~\ref{sec:noise}. To quantify parameter stability across runs, we adopt a metric based on the signal-to-noise ratio. Such variability analysis is related to identifiability issues in nonlinear system identification \cite{ljung1994identifiability}.
For a given learned parameter $\theta$ estimated across $K$ independent training runs, we define
\[
\mathrm{SNR}(\theta) = \frac{\mu_\theta}{\sigma_\theta},
\]
where $\mu_\theta$ and $\sigma_\theta$ denote the sample mean and standard deviation of $\theta$ across runs, respectively. Large values indicate that the parameter is consistently estimated and well concentrated around its mean, while small values signal that training variability is comparable to or exceeds the identified level. The following thresholds are adopted for interpretation: $\mathrm{SNR} < 1$ indicates that variability dominates and the parameter is not stably identified; $1 \leq \mathrm{SNR} < 2$ denotes marginal stability; $\mathrm{SNR} \geq 2$ corresponds to a structurally stable estimate; and $\mathrm{SNR} \geq 3$ indicates strong concentration around the mean.

The performance metrics across the ten runs exhibit very high reproducibility. Training loss, validation loss, test RMSE, and test PCC all present SNR values well above $30$, confirming that the trajectory-based loss leads to consistent convergence regardless of the specific noise realization in the dataset.

At the parameter level, the thermal graph comprises $37$ learnable parameters, collecting nodal capacitances $\gamma_i$ and edge conductances $\delta_{ij}$. Of these, $31$ exceed the stability threshold of $\mathrm{SNR} \geq 3$. Six parameters fall below this threshold, of which two present $\mathrm{SNR} < 2$, with a minimum observed value of $1.1$.

\smallskip

The larger variability observed for a small subset of parameters is consistent with limited information content in the available excitation/measurement setup. As a consequence, multiple parameter values can yield nearly indistinguishable trajectories on the observed components, leading to larger run-to-run variability without affecting predictive accuracy. This behavior reflects intrinsic practical identifiability limitations of the inverse problem, rather than instability of the training procedure.
\section{Conclusions}
\label{sec:6}

In this work, we proposed a physics-constrained graph-based identification framework for Lumped Parameter Thermal Models (LPTMs), designed to bridge the gap between physics-based reduced-order models and purely data-driven approaches for spacecraft thermal systems. The methodology formulates thermal dynamics within the formalism of graph-based Neural ODEs, while preserving a physically interpretable message-passing structure and without relying on black-box neural parametrizations, enabling continuous-time parameter identification directly from temperature measurements and known thermal inputs.

A key feature of the proposed approach is the explicit enforcement of the physical structure of the thermal network. The model is constructed so that all admissible parameters satisfy by construction fundamental properties of conductive heat transfer, including positivity of nodal coefficients, symmetry of interactions, and dissipativity of the resulting dynamics. These properties are enforced by construction through a bounded parameterization, without requiring additional regularization terms. The inverse calibration problem is formulated as a trajectory-matching optimization over a graph dynamical system, and solved within a Neural ODE framework using adjoint-based gradient computation. To ensure numerical stability and computational efficiency, the training pipeline combines bounded parameter embeddings, trajectory batching through block-diagonal graph replication, and controlled forward evolution. This design enables robust training even in the presence of long time horizons and stiff thermal dynamics.

The proposed framework was validated on synthetic datasets generated from high-fidelity FEM simulations under progressively complex forcing conditions. The results demonstrate accurate long-horizon predictions, with high correlation between predicted and reference trajectories and stable autoregressive behavior over multiple orbital cycles. Sensitivity analyses show limited degradation under measurement noise and highlight the role of the training window length in balancing optimization difficulty and generalization. Repeated training runs confirm that most parameters are consistently identified, while residual variability is confined to weakly excited thermal interactions, reflecting intrinsic identifiability limitations rather than instability of the method.

From a modeling perspective, the framework provides a modular, interpretable, and physically grounded approach to the calibration of reduced-order thermal models. The identified parameters retain a clear physical meaning as effective thermal properties, while the graph-based Neural ODE formulation enables flexible integration of data and prior structural knowledge. This combination allows the construction of computationally efficient models suitable for predictive tasks in operational environments.

In the context of spacecraft thermal Digital Twins, the proposed methodology provides a promising basis for continuous calibration from telemetry data, supporting real-time monitoring, anomaly detection, and predictive thermal control. The resulting models can be used to assess deviations from nominal behavior and to evaluate alternative control strategies under varying operational scenarios. The proposed identification framework aligns with current digital‑twin initiatives for spacecraft thermal management \cite{glaessgen2012digital,scigliano2024heatpipe}.

Validation of the presented approach on real telemetry data remains an important direction for future work,together with its integration into onboard real-time temperature control systems.

\section*{Acknowledgments}

This work has been carried out in the framework of the Hybrid Digital Twins for Satellites (HDTS) Work Package of the ASTRAL project, funded by the Research National Center in High Performance Computing, Big Data and Quantum Computing through an Innovation Grant assigned to Thales Alenia Space Italia (TASI). The authors gratefully acknowledge the support of TASI technical staff whose domain expertise was crucial in defining the reference use case, while their continuous support provided essential system-level insights into the spacecraft's thermal subsystem, ensuring the industrial relevance of the proposed framework. SP, NP, FR, LS  are members of INdAM-GNCS group. The present research is part of the activities of “Dipartimento di Eccellenza 2023-2027”. SP, FR, LS acknowledge support from project FIS, MUR, Italy 2025-2028, Project code: FIS-2023-02228, CUP: D53C24005440001, ``SYNERGIZE: Synergizing Numerical Methods and Machine Learning for a new generation of computational models''.

\newpage



\bibliographystyle{elsarticle-num} 
\bibliography{references}          

\end{document}